%
%
%




%

\count100= 1
\count101= 10

\magnification\magstep1
\hsize 4.43in
\vsize 7.5in



\font\seventeenrm=cmr17
\font\twelverm=cmr12
\font\ninerm=cmr9
\font\sevenrm=cmr7

\font\fiverm=cmr5
\font\ninei= cmmi9
\font\seveni= cmmi7
\font\ninesy= cmsy9
\font\sevensy= cmsy9

\font\ninebf = cmbx9
\font\nineit = cmti9
\font\ninesl = cmsl9
\font\nineex = cmex9

\font\efont=cmti10

\def\figfont{\ninerm%
    \textfont0 = \ninerm
    \textfont1 = \ninei
    \textfont2 = \ninesy
    \textfont3 = \nineex
    \scriptfont0 = \sevenrm
    \scriptfont1 = \seveni
    \scriptfont2 = \sevensy
    \scriptscriptfont0 = \fiverm
    \scriptscriptfont1 = \fivei
    \scriptscriptfont2 = \fivesy
    \let \rm = \ninerm
    \let \sl = \ninesl
    \let \bf = \ninebf
    \let \it = \nineit
    \baselineskip 9pt}


\def\CC{{\rm C\kern-.18cm\vrule width.6pt height 6pt depth-.2pt
\kern.18cm}}

\def\NN{{\mathop{{\rm I}\kern-.2em{\rm N}}\nolimits}}

\def\PP{{\mathop{{\rm I}\kern-.2em{\rm P}}\nolimits}}

\def\RR{{\mathop{{\rm I}\kern-.2em{\rm R}}\nolimits}}

\def\RRt{{\seventeenrm I}\kern-.2em{\seventeenrm R}}


\def\ZZ{{\mathop{{\rm Z}\kern-.28em{\rm Z}}\nolimits}}


\def\bfm#1{{\dimen0=.01em\dimen1=.009em\makebold{$#1$}}}



\def\makebold#1{\mathord{\setbox0=\hbox{#1}%
       \copy0\kern-\wd0%
       \raise\dimen1\copy0\kern-\wd0%
       {\advance\dimen1 by \dimen1\raise\dimen1\copy0}\kern-\wd0%
       \kern\dimen0\raise\dimen1\copy0\kern-\wd0%
       {\advance\dimen1 by \dimen1\raise\dimen1\copy0}\kern-\wd0%
       \kern\dimen0\raise\dimen1\copy0\kern-\wd0%
       {\advance\dimen1 by \dimen1\raise\dimen1\copy0}\kern-\wd0%
       \kern\dimen0\raise\dimen1\copy0\kern-\wd0%
       \kern\dimen0\box0}}




\def \eword#1{{\efont #1}}

\def\eg{{\it e.g.}}%

\def\frac#1#2{{#1 \over #2}}


\def\ie{{\it i.e.}}

\def\ms{\medskip}

\def\noin{\noindent}









\def\eop{\makeblanksquare6{.4}\ms}


\def\nopf{\medskip}

\def\pf{\noindent{\bf Proof: }}

\def\makeblanksquare#1#2{
\dimen0=#1pt\advance\dimen0 by -#2pt
      \vrule height#1pt width#2pt depth0pt\kern-#2pt
      \vrule height#1pt width#1pt depth-\dimen0 \kern-#1pt
      \vrule height#2pt width#1pt depth0pt \kern-#2pt
      \vrule height#1pt width#2pt depth0pt
}


\def\abstract#1{\bigskip\bigskip\medskip%
 {\narrower\noindent\figfont{\bf Abstract.~}#1\bigskip}\medskip}%

\def\author#1{\bigskip\bigskip\centerline{\twelverm #1}}

\def\footnoterule{\kern -3pt \hrule width 0truein \kern 2.6pt}
\def\leftheadline{\ifnum\pageno=\count100 \hfill%
  \else\rm\folio\hfil\it\shortauthor\fi}
\def\rightheadline{\ifnum\pageno=\count100 \hfill%
  \else\it\shorttitle\hfil\rm\folio\fi}
\def\title#1{\centerline {\seventeenrm #1}}
\def\titwo#1{\medskip \centerline {\seventeenrm #1}}

\def\titexp#1#2{\hbox{{\seventeenrm #1} \kern-.25em%
  \raise .90ex \hbox{\twelverm #2}}\/}
\def\titsub#1#2{\hbox{{\seventeenrm #1} \kern-.25em%
  \lower .60ex \hbox{\twelverm #2}}\/}

\nopagenumbers
\headline{\ifodd\pageno\rightheadline \else\leftheadline\fi}
\footline{\hfil}
\null\vskip 18pt
\centerline{}
\pageno=\count100
\count102=\count100
\advance\count102 by -1
\advance\count102 by \count101


\def\copyright{\hbox{{\twelverm o}\kern-.61em\raise .46ex\hbox{\fiverm c}}}



\def\sect#1{\goodbreak\bigskip\smallskip\centerline{\bf\S #1}\medskip
    \noindent\ignorespaces}






\newbox\maboite

\def\boxit#1#2{\setbox\maboite=\hbox{\kern#1{#2}\kern#1}%
    \dimen1=\ht\maboite \advance\dimen1 by #1 \dimen2=\dp\maboite
\advance\dimen2 by #1%
    \setbox\maboite=\hbox{\vrule height\dimen1
depth\dimen2\box\maboite\vrule}%
    \setbox\maboite=\vbox{\hrule\box\maboite\hrule}%
    \advance\dimen1 by .4truept \ht\maboite=\dimen1%
    \advance\dimen2 by .4truept \dp\maboite=\dimen2 \box\maboite\relax}


\def\Acknowledgments{\goodbreak\bigskip\noindent{\bf
   Acknowledgments.\ }}

\def\Remark#1.{\goodbreak\medskip\noin {\bf Remark#1.}}
\def\Example#1.{\goodbreak\medskip\noin {\bf Example#1.}}


\def\ref{\smallskip\global\advance\refnum by 1 \item{\the\refnum.}}
\newcount\refnum \refnum = 0

\def\References{\goodbreak\bigskip\centerline{\bf References}\bigskip
   \frenchspacing}





\title{Convex Multivariate Approximation}
\titwo{by Algebras of Continuous Functions}
\author{Andriy V. Bondarenko and Andriy V. Prymak}

\def\shorttitle{Convex Multivariate Approximation}
\def\shortauthor{A. V. Bondarenko and A. V. Prymak}


\def\bp{{SG-}}

\def\norm#1#2{\left\|#1\right\|_{#2}}



\abstract {We obtain an analog of Shvedov theorem for
convex multivariate approximation by algebras of continuous
functions.}


\sect{1. Introduction}
Let $K\subset \RR^d$, $d\in\NN$, be a compact set and $C(K)$ be
the space of continuous functions $f$ on $K$ with the uniform norm
$$
\norm{f}{C(K)}=\max_{{\bfm x}\in K}|f({\bfm x})|.
$$
For $X\subset C(K)$, denote by $\overline{X}$ the closure of $X$
with respect to this norm. Recall, a continuous function $f$ on a
convex set $K_*\subset\RR^d$ is called a convex function on $K_*$,
if for all ${\bfm x_1}, {\bfm x_2}\in K_*$ the inequality
$$
\frac {f({\bfm x_1})+f({\bfm x_2})}2\ge f\left (\frac{{\bfm
x_1}+{\bfm x_2}}2\right) \eqno (1)
$$
holds. Then Shvedov proved the following analog of Weierstrass
theorem for convex multivariate approximation.
\proclaim Theorem S[3].
Let $K\subset\RR^d$ be a convex compact set. If $f$ is a convex
function on $K$, then for each $\epsilon>0$ there exists an
algebraic polynomial $p$ of $d$ variables, such that $p$ is convex
on $\RR^d$, and $\norm{p-f}{C(K)}<\epsilon$. \nopf

Applying Theorem~S we obtain its analog for convex multivariate
approximation by algebras of continuous functions.
 To formulate our result, we introduce a definition.
First we recall, that a set $A$ of functions $f\in C(K)$ is called
a real algebra with unity, if $A$ is a real linear space,
satisfying: if $f\in A$ and $g\in A$, then $fg\in A$, and $f_0\in
A$, where $f_0({\bfm x})\equiv 1$.
\proclaim Definition.
Let $\Delta\subseteq C(K)$. A set $D\subseteq\overline{\Delta}$ is
said to be the \bp set (shape generating) of $\Delta$ if
\eword{for all} real algebras with unity $A\subseteq C(K)$ such that
$D\subseteq\overline{A\cap\Delta}$, we have
$\Delta\subseteq\overline{A\cap\Delta}$. \nopf

In other words, if $D$ is the \bp set of some set $\Delta$ of
constrains,
 then in order to determine whether the elements of an
algebra can provide constrained approximation to any function from
$\Delta$, it is sufficient to verify the possibility of
constrained approximation by the elements of the algebra to any
function from $D$. If $D$ is ``much smaller'' than $\Delta$, this
provides an efficient condition on the possibility of constrained
approximation by the elements of the algebra.

 Recall, that by
classical Stone-Weierstrass theorem (see say [1, p.11]), if
$K\subset\RR^d$ is a compact set and $A\subseteq C(K)$ is an
algebra with unity, then $\overline{A}=C(K)$ if and only if $A$
separates the points of $K$, \ie, for each pair
$\bfm{x},\bfm{y}\in K$, such that $\bfm{x}\ne\bfm{y}$, there is a
function $a\in A$ satisfying $a(\bfm{x})\ne a(\bfm{y})$. It is not
hard to see, that this result can be stated in the terms of \bp
sets: {\sl For any compact set $K\subset\RR^d$, the set
$\{l_1,\ldots,l_d\}$ is an
\bp set of $\Delta=C(K)$, where
$$
l_j(x_1,\dots,x_d)=x_j, \quad j=1,\dots,d.\eqno(2)
$$ }

We find  \eword{finite} \bp sets for the convex multivariate
approximation, and prove that the cardinality of such \bp sets can
not be reduced. Our main result is
\proclaim Theorem 1. Let
$K\subset \RR^d$ be a convex compact set, and $\Delta$ be the set
of convex and continuous functions on $K$. Then

{\sl a) the set
$$
D:=\{l_1,\dots,l_{d+1}\},\eqno(3)
$$
where $l_j$, $j=1,\ldots,d$, are given by~(2), and
$$
l_{d+1}(x_1,\dots,x_d)=-(x_1+\dots+x_d),
$$
is an \bp set of $\Delta$,

 b) if $K$ has non-empty interior, then
the minimal cardinality of \bp set of $\Delta$ is equal to
$d+1$.\sl}\nopf

This theorem consists of two parts: positive result and
counterexample. The proof of the positive part is constructive. In
fact, we can choose arbitrary $d+1$ linear functions
$l_1,\dots,l_{d+1}$, such that each linear function of $d$
variables can be represented in the form
$\lambda_0+\lambda_1l_1+\dots+\lambda_{d+1}l_{d+1}$, where
$\lambda_1,\dots,\lambda_{d+1}\ge0$. For the counterexample we use
the well-known Borsuk antipodality
\proclaim Theorem B[2, Appendix 1].
Let $S^{d-1}$ be the unit sphere of $\RR^d$. If
$f:S^{d-1}\to\RR^{d-1}$ is a continuous and odd mapping, then
there exists a point ${\bfm x}\in S^{d-1}$, such that $f({\bfm
x})=0$.
\nopf

We prove Theorem~1 in Section~2 and Section~3. Examples are given
in Section~4.


\sect{2. Proof of Theorem~1a)}
Recall, that $\Delta$ is the set of all convex and continuous
functions on some convex compact $K\subset \RR^d$. In our proof we
often use the obvious fact that $\Delta$ is a convex cone, \ie, if
$f,g\in\Delta$ and $\alpha,\beta\ge0$, then $\alpha f+\beta
g\in\Delta$. We have to show that the set $D$ defined by~(3) is a
\bp set of $\Delta$. To this end first we prove \proclaim Lemma 1.
Let $K_*\subset\RR^d$ be a convex set, $g:K_*\to\RR$ be a convex
function on $K_*$ and non-decreasing in each variable,
$f_j:K\to\RR$, $j=1,\dots,d$, be convex functions, such that for
all $\bfm{x}\in K$ we have $(f_1(\bfm{x}),\dots,f_d(\bfm{x}))\in
K_*$. Then $g(f_1,\dots,f_d)$ is convex function on $K$.

\pf For any $\bfm{x_1},\bfm{x_2}\in K$ we have
$$
\eqalign{\frac12
(g(f_1(\bfm{x_1}),&\dots,f_d(\bfm{x_1}))+
g(f_1(\bfm{x_2}),\dots,f_d(\bfm{x_2}))) \cr &\ge
g\left(\frac{f_1(\bfm{x_1})+f_1(\bfm{x_2})}2,\dots,
\frac{f_d(\bfm{x_1})+f_d(\bfm{x_2})}2\right)
\cr &\ge g\left(f_1\left(\frac{\bfm{x_1}+\bfm{x_2}}2\right),
\dots,f_d\left(\frac{\bfm{x_1}+\bfm{x_2}}2\right)\right),
\cr }
$$
where in the first inequality we use the convexity of $g$ and in
the second inequality we use the convexity of $f_j$, $j=1,\dots,d$
and the fact that $g$ is non-decreasing in each variable.
\eop

Fix an algebra with unity $A\subset C(K)$, satisfying
$D\subset\overline{A\cap\Delta}$, and put
$$
K_*:=\{\bfm{x}\in\RR^d|\,\exists \bfm{y}\in K :
\rho(\bfm{x},\bfm{y}):=\Bigl(\sum_{j=1}^d(x_j-y_j)^2\Bigr)^{1/2}\le1\}.
$$
Since $K$ is a convex compact set, the set $K_*$ is a convex
compact set containing $K$. By Theorem~S, it is sufficient to
check that all convex algebraic polynomials on $\RR^d$ belong to
$\overline{A\cap\Delta}$. Indeed, let $p$ be an arbitrary convex
algebraic polynomial on $\RR^d$. We set
$$
l(\bfm{x}):=\sum_{j=1}^d\norm{\frac{\partial p}{\partial x_j}}
{C(K_*)}x_j,
$$
where $\bfm{x}=(x_1,\dots,x_d)\in\RR^d$. Since $-l$ can be
represented as a linear combination of functions $l_j$,
$j=1,\ldots ,d+1$, with positive coefficients, then
$-l\in\overline{A\cap\Delta}$. So, to complete the proof we have
to show that
$$
g:=p+l\in\overline{A\cap\Delta}.
$$
We remark, that $g$ is constructed in this way to be
non-decreasing in each variable on $K_*$. Given $\epsilon>0$, let
$\delta\in (0,1)$ be such, that $\|g(\bfm{x_1})-
g(\bfm{x_2})\|_{C(K_*)}<\epsilon$, if
$\rho(\bfm{x_1},\bfm{x_2})<\delta$, and $\bfm{x_1},\bfm{x_2}\in
K_*$ (here we use the fact, that any continuous function on a
compact set is uniformly continuous). Since
$D\subset\overline{A\cap\Delta}$, then for each $j=1,\ldots,d$
there is a function $h_j\in A\cap\Delta$, such that
$\|h_j-l_j\|_{C(K)}<\delta/\sqrt{d}$. Put $\overline{\bfm
{h}}(\bfm{x}):= (h_1(\bfm{x}),\ldots, h_d(\bfm{x}))$, and
$\overline{\bfm{l}}(\bfm{x}):=(l_1(\bfm{x}),\ldots,
l_d(\bfm{x}))\equiv \bfm{x}$. Then $\rho(\overline{\bfm
{h}}(\bfm{x}),\bfm{x})=\rho(\overline{\bfm
{h}}(\bfm{x}),\overline{\bfm{l}}(\bfm{x}))<\delta<1$ for all
$\bfm{x}\in K$, in particular $\overline{\bfm {h}}(\bfm{x})\in
K_*$, for $\bfm{x}\in K$. Therefore, with notation
$g\circ\overline{\bfm {h}}=g(h_1,\ldots,h_d)$, we have
$$
\|g\circ\overline{\bfm {h}}-g\|_{C(K)}<\epsilon.
$$
Evidently, $g$ is a polynomial, hence
$$
g\circ\overline{\bfm {h}}\in A.
$$
Finally, definitions of $l$ and $g$ imply, that $g$ is convex on
$K_*$ and non-decreasing in each variable. Therefore, Lemma~1
yields
$$
g\circ\overline{\bfm {h}}\in\Delta.
$$
For a given sequence $\epsilon_n\to0$, $n\to\infty$, construct a
sequence $\overline{\bfm {h}}_n$ such that $\|g\circ\overline{\bfm
{h}}_n-g\|_{C(K)}<\epsilon_n$. We have $g\circ\overline{\bfm
{h}}_n\to g$, as $n\to\infty$, and so
$g\in\overline{A\cap\Delta}$. Part~a) of Theorem~1 is proven.

\sect{3.Proof of Theorem~1b)}
For each $f_1,\dots,f_d\in\Delta$ we construct an example of
algebra with unity $A:=A(f_1,\dots,f_d)$ such that
$\overline{A\cap\Delta}\supset\{f_1,\ldots,f_d\}$ and
$\overline{A\cap\Delta}\ne\Delta$.

Without loss of generality one may assume, that $0$ is an interior
point of $K$. For each $\bfm{y}\in S^{d-1}$ set
$I_\bfm{y}:=\{t\in\RR : t\bfm{y}\in K\}$. Since $0$ is an interior
point of $K$, $I_\bfm{y}$ is a closed interval with non-empty
interior. Moreover, each point $t\bfm y$, $t\in{\rm Int}~I_{\bfm
y}$, is an interior point of $K$ (here and below ${\rm Int}~M$
denotes the interior of $M$). We need
\proclaim Lemma 2. Suppose $h\in C(K)$, a sequence ${\bfm
{y_0}},{\bfm {y_1}},\ldots\in S^{d-1}$ satisfies ${\bfm
{y_k}}\to{\bfm {y_0}}$, as $k\to\infty$, and numbers $m_k$,
$k=1,2,\dots$, are such that
$$
\min_{t\in I_{{\bfm {y_k}}}}h(t{\bfm {y_k}})=h(m_k{\bfm {y_k}}),
$$
and $m_k\to m_0$, as $k\to\infty$. Then $m_0{\bfm {y_0}}\in K$,
and
$$
\min_{t\in I_{{\bfm {y_0}}}}h(t{\bfm {y_0}})=h(m_0{\bfm {y_0}}).
$$

\pf
We have $m_k{\bfm {y_k}}\to m_0{\bfm {y_0}}$, as $k\to\infty$, so,
due to the compactness of $K$, we get $m_0{\bfm {y_0}}\in K$. Let
$t_*\in I_{{\bfm {y_0}}}$ be fixed. We claim, that one can choose
a sequence of numbers $t_1,t_2,\ldots$, satisfying $t_k\in
I_{{\bfm {y_k}}}$, $k=1,2,\dots$, and $t_k\to t_*$, as
$k\to\infty$. Indeed, it is sufficient to have $I_{{\bfm
{y_k}}}\to I_{{\bfm {y_0}}}$, as $k\to\infty$, in the sense that
the endpoints of $I_{{\bfm {y_k}}}$ tend to the corresponding
endpoints of $I_{{\bfm {y_0}}}$ (recall that all $I_{\bfm{y_j}}$,
$j=0,1,2,\ldots$ are closed segments). This follows from the fact
that the boundary of convex compact $K$ is a closed set, and
$t\bfm{y_j}\in {\rm Int}~K$, for $t\in{\rm Int}~I_{\bfm{y_j}}$,
$j=0,1,2,\ldots$. So, there are $t_k\in I_{{\bfm {y_k}}}$,
$k=1,2,\dots$, such that $t_k\to t_*$. We have $h(t_k{\bfm
{y_k}})\ge h(m_k{\bfm {y_k}})$, $k=1,2,\ldots$, and, therefore,
$h(t_*{\bfm {y_0}})\ge h(m_0{\bfm {y_0}})$.
\eop

Now we prove
\proclaim Proposition 1. For
arbitrary $d$ functions $f_1,\dots,f_d\in\Delta$ there exists a
point ${\bfm y}\in S^{d-1}$ and a number $m\in I_{{\bfm y}}$, such
that
$$
\min_{t\in I_{{\bfm{y}}}}f_j(t\bfm{y})=f_j(m\bfm{y}), \quad
{\sl for~all~} j=1,\dots,d.
$$

\pf For $d=1$ this is trivial. Let $d\ge2$. Suppose
$f_1,\dots,f_d$ are strictly convex on $K$, that is
 inequality (1)
holds with ``$>$'' instead of ``$\ge$'', for all
$\bfm{x_1}\ne\bfm{x_2}$. For any strictly convex function $h$ and
each $\bfm{y}\in S^{d-1}$ we denote by $m:=m(h,\bfm{y})$ the
unique point, satisfying
$$
h(m\bfm{y})=\min_{t\in I_\bfm{y}}h(t\bfm{y}),
$$
Clearly, for each $h$ $m(h,\cdot)$ is an odd mapping from
$S^{d-1}$ to $\RR$. We claim that $m(h,\cdot)$ is continuous.
Indeed, otherwise, there exists a sequence ${\bfm {y_0}},{\bfm
{y_1}},\ldots\in S^{d-1}$, such that ${\bfm {y_k}}\to{\bfm
{y_0}}$, as $k\to\infty$, and for $m_k:=m(h,{\bfm {y_k}})$,
$k=1,2,\dots$, $m_*:=m(h,{\bfm {y_0}})$, we have $m_k\not\to m_*$,
as $k\to\infty$. The sequence ${m_1},{m_2},\ldots$, is bounded
(by,
\eg, the diameter of $K$), hence, there is a
subsequence ${m_{k(1)}},{m_{k(2)}},\ldots$, convergent to some
$m_0$, $m_0\ne m_*$. By Lemma~2, $m_0{\bfm {y_0}}\in I_{\bfm
{y_0}}$, and $h(m_*{\bfm {y_0}})\ge h(m_0{\bfm {y_0}})$, we get
$h(m_*{\bfm {y_0}})=h(m_0{\bfm {y_0}})$, which contradicts the
strict convexity of $h$.

For $\bfm{y}\in S^{d-1}$ we define
$$
\eqalign{
g(&\bfm{y}):= \cr &(m(f_1,\bfm{y})-m(f_d,\bfm{y}),
m(f_2,\bfm{y})-m(f_d,\bfm{y}),\dots,
m(f_{d-1},\bfm{y})-m(f_d,\bfm{y})).\cr}
$$
The mapping $g$ from $S^{d-1}$ to $\RR^{d-1}$ is continuous and
odd, hence, by Borsuk antipodality theorem, there is a point
$\bfm{y'}\in S^{d-1}$ such that $g(\bfm{y'})=0$. We obtain
$m(f_1,\bfm{y'})=m(f_2,\bfm{y'})=\dots=m(f_d,\bfm{y'})$, and the
proposition is proven for strictly convex functions.

For the general case we approximate each $f_j$ by a sequence of
strictly convex functions
$$
f_{j,n}(\bfm{x}):=f_j(\bfm{x})+\frac1n\rho^2(\bfm{x},\bfm{z}),
$$
where $\bfm{z}\in\RR^d$ is an arbitrary fixed point. For each
$n\in\NN$ and functions $f_{1,n},\ldots, f_{d,n} \in \Delta $ we
get the corresponding ${\bfm {y}_n}$ and the number $m_n$. Since
$S^{d-1}$ is a compact set and $m_n$ are bounded, we can choose a
convergent subsequence $({m_{n({k})}},{\bfm{y}_{n({k})}})\to
(m,\bfm{y})$, as $k\to\infty$. The relation
$$
\eqalign{
\min_{t\in I_{\bfm y}}(f_j(t \bfm{y}_{n(k)})&+\rho^2(t{\bfm y}_{n(k)},\bfm
z)/n(k))\cr&= f_j(m_{n(k)} \bfm{y}_{n(k)})+\rho^2(m_{n(k)}{\bfm
y}_{n(k)},\bfm z)/n(k)\cr }
$$
implies
$$
f_j(m_{n(k)}{\bfm{y}}_{n(k)})+\rho^2(m_{n(k)}{\bfm y}_{n(k)},\bfm
z)/n(k)\le f_j(t{\bfm{y}}_{n(k)})+\rho^2(t{\bfm y}_{n(k)},\bfm
z)/n(k)
$$
for all $t \in I_{{\bfm y}_{n(k)}}$. Taking limits when
$k\to\infty$ and taking into account that any $t\in I_{\bfm{y}}$
is a limit of $t\in I_{{\bfm y}_{n(k)}}$, we obtain
$$
f_j(m{\bfm y})\le f_j(t\bfm y)
$$
and
$$
\min_{t\in I_{{\bfm {y}}}}f_j(t{\bfm {y}})=f_j(m{\bfm {y}}),
\quad {\rm for~all~} j=1,\dots,d.
$$
Proposition 1 is proven. \eop

Fix arbitrary $d$ functions $f_1,\ldots,f_d\in \Delta$. Applying
Proposition~1 for them, we obtain a point $\bfm{y}$  and a number
$m$. Set $I:=\{t\bfm{y}:t\in I_{\bfm{y}}\}$, and $\bfm{x'}:=m\bfm
y$. Denote $A(\bfm{x'},I)$ the set of all functions $f$ continuous
on $K$ with the property: there exists a closed interval
$J\subseteq I$ with non-empty interior, satisfying $\bfm{x'}\in J$
and $f(\bfm{x})=f(\bfm{x'})$ for all $\bfm{x}\in J$. Clearly,
$A(\bfm{x'},I)$ is an algebra with unity. With the above
notations, we prove
\proclaim Proposition 2.
Put $\Delta':=\{f\in\Delta:\min_{\bfm{x}\in
I}f(\bfm{x})=f(\bfm{x'})\}$. We have
$\overline{A(\bfm{x'},I)\cap\Delta}=\Delta'$.

\pf  The properties of convex functions and the definitions of $A$
and $\Delta'$ easily yield
$A(\bfm{x'},I)\cap\Delta\subset\Delta'$, hence
$\overline{A(\bfm{x'},I)\cap\Delta}\subset\overline{\Delta'}=\Delta'$.
Let us prove that
$\overline{A(\bfm{x'},I)\cap\Delta}\supset\Delta'$. Let
$f\in\Delta'$ be a convex function. Given $\epsilon>0$, put
$$
\eqalign{
L_\epsilon:=\{(x_1,&\ldots,x_{d+1})\in\RR^{d+1}:\cr
&(x_1,\ldots,x_d)=\lambda\bfm{y},\lambda\in\RR,
x_{d+1}=f(\bfm{x'})-\epsilon\},\cr}
$$
and
$$
\eqalign{F:=\{(x_1,&\ldots,x_{d+1})\in\RR^{d+1}:\cr &(x_1,\ldots,x_d)\in K,
f(x_1,\ldots,x_d)\le x_{d+1}\le M\},\cr}
$$
where
$$
M:=\max_{\bfm{x}\in K}f(\bfm{x}).
$$
Clearly, $L_\epsilon$ is a straight line in $\RR^{d+1}$. $F$ is
the ``truncated'' epigraph of $f$, and $F$ is a convex set which
is the intersection of the epigraph of $f$ and the half-space
$\{(x_1,\ldots,x_{d+1}):x_{d+1}\le M\}$. So, both $L_\epsilon$ and
$F$ are closed and convex sets, moreover $F$ is a compact set.
Since $f\in\Delta'$, we have $L_\epsilon\cap F=\emptyset$. Hence,
there is a hyperplane $s_*$, separating $L_\epsilon$ and $F$. In
other words, this means that for some reals $a_0,\ldots,a_{d+1}$,
such that $|a_1|+\ldots+ |a_{d+1}|\ne0$ and $a_{d+1}\ge0$, we
have
$$
s_*=\{(x_1,\ldots,x_{d+1})\in\RR^{d+1}:a_1x_1+\ldots+a_{d+1}x_{d+1}=a_0\},
$$
and the half-spaces
$$
s_*^+=\{(x_1,\ldots,x_{d+1})\in\RR^{d+1}:a_1x_1+\ldots+a_{d+1}x_{d+1}>
a_0\},
$$
and
$$
s_*^-=\{(x_1,\ldots,x_{d+1})\in\RR^{d+1}:a_1x_1+\ldots+a_{d+1}x_{d+1}<
a_0\},
$$
satisfy
$$
s_*^+\supset F,\quad s_*^+\cap L_\epsilon=\emptyset,\quad
s_*^-\supset L_\epsilon,\quad s_*^-\cap F=\emptyset.\eqno(4)
$$
Let $\bfm{x'}=(x_1',\ldots,x_d')$. Consider the line
$$
L:=\{(x_1,\ldots,x_{d+1})\in\RR^{d+1}:x_1=x_1',\ldots,x_d=x_d',x_{d+1}\in\RR\}.
$$
If $a_{d+1}=0$, then either $L\subset s_*^+$, or $L\subset s_*^-$.
But both $(x_1',\ldots,x_d',f(\bfm{x'}))$ and
$(x_1',\ldots,x_d',f(\bfm{x'})-\epsilon)$ belong to $L$, and, in
the same time, we have $(x_1',\ldots,x_d',f(\bfm{x'}))\in F$,
$(x_1',\ldots,x_d',f(\bfm{x'})-\epsilon)\in L_\epsilon$, which
contradicts~(4). So, $a_{d+1}>0$. This implies, that $s_*$ is a
graph of some linear function $\tilde{s}:\RR^d\to\RR$,
$$
\tilde{s}(x_1,\ldots,x_d)=a_{d+1}^{-1}(a_0-a_1x_1-\ldots-a_dx_d),
\quad (x_1,\ldots,x_d)\in\RR^{d},
$$
satisfying $\tilde{s}(\bfm{x})<f(\bfm{x})$, $\bfm{x}\in K$, and
$\tilde{s}(\lambda\bfm{y})>f(\bfm{x'})-\epsilon$, $\lambda\in\RR$.
The latter provides $\tilde{s}(\lambda\bfm{y})=const$,
$\lambda\in\RR$, since $\tilde{s}$, restricted to the line
$\{\lambda\bfm{y}:\lambda\in\RR\}$, is a linear function. Put
$s(\bfm{x}):=\tilde{s}(\bfm{x})+f(\bfm{x'})-\epsilon-\tilde{s}(\bfm{x'})$,
$\bfm{x}\in K$. We have, that $s$ is a hyperplane, satisfying
$f(\bfm x)> s(\bfm x)$, $\bfm x\in K$, $s({\bfm x'})=f({\bfm
x'})-\epsilon$, $s(\bfm x)=s({\bfm x'})$, $\bfm{x}\in I$.
 Put $J_{\epsilon}:=\{{\bfm x}\in I|
f(\bfm{x})\le f(\bfm{x'})+\epsilon\}$, and for ${\bfm x}\in K$
$$
g_{\epsilon}({\bfm x}):=\max\{f(\bfm x)-2\epsilon,s(\bfm x)\}.
$$
Clearly, $g_\epsilon\in\Delta$. We have
$g_\epsilon(\bfm{x})=f(\bfm{x'})-\epsilon$, $\bfm{x}\in
J_\epsilon$. Indeed, if $\bfm{x}\in J_\epsilon$, then
$f(\bfm{x})-2\epsilon\le f({\bfm {x'}})-\epsilon=s({\bfm
{x'}})=s(\bfm{x})$, and
$g_\epsilon(\bfm{x})=s(\bfm{x})=s(\bfm{x'})$. So, $g_\epsilon\in
A(\bfm{x'},I)$.
 The definition of
$g_\epsilon$ implies that for $\bfm{x}\in K$ we have
$g_\epsilon(\bfm{x})\ge f(\bfm{x})-2\epsilon$, and
$g_\epsilon(\bfm{x})=\max\{f(\bfm x)-2\epsilon,s(\bfm x)\}<
\max\{f(\bfm x)-2\epsilon,f(\bfm x)\}=f(\bfm{x})$, so $f(\bfm{x})>
g_{\epsilon}(\bfm{x}) \ge f(\bfm{x})-2\epsilon$. This implies $
f\in\overline{A(\bfm{x'},I)\cap\Delta}. $
\eop
To conclude the proof of Theorem~1b), we remark that
$f_1,\ldots,f_d\in\Delta'$, and $\Delta'\ne\Delta$.


\sect{4.Examples}
First we give an easy result, allowing to reveal a vast amount of
univariate algebras, providing convex approximation.
\proclaim Proposition~3. Let $A\subset C[a,b]$ be an arbitrary algebra of
univariate functions, containing a twice continuously
differentiable function $h\in C^{(2)}[a,b]$, satisfying $h'(x)>0$,
for all $x\in[a,b]$. Let $\Delta$ be the set of all convex and
continuous functions on $[a,b]$. Then $\Delta=\overline{\Delta\cap
A}$.

\pf It is sufficient to prove that
$p\in\overline{\Delta\cap A}$ for all polynomials $p\in\Delta$. To
this end, fix a polynomial $p\in\Delta$ and $\delta>0$, and put
$p_\delta(x):=p(x)+\delta x^2/2$, $x\in[a,b]$. We have
$p_\delta''(x)\ge\delta$, $x\in[a,b]$.

Since $h$ is strictly monotone and continuous, there is the
inverse function $h^{-1}:[a_1,b_1]\to[a,b]$, where $a_1=h(a)$,
$b_1=h(b)$. Moreover, since $h'>0$, the condition $h\in
C^{(2)}[a,b]$ implies $h^{-1}\in C^{(2)}[a_1,b_1]$.

Put $q:=p\circ h^{-1}$, \ie, $q(x):=p_\delta(h^{-1}(x))$,
$x\in[a_1,b_1]$. We have that $q$ is a twice continuously
differentiable function, therefore, there is a sequence of
polynomials $p_n$, $n=1,2,\ldots$, satisfying
$$
\norm{p_n^{(j)}-q^{(j)}}{C[a_1,b_1]}\to0,\quad {\rm as~~}
n\to\infty, \quad j=0,1,2. \eqno(5)
$$
It follows, that $(p_n\circ h)''=(p_n''\circ h) h'^2+(p_n'\circ
h)h''$ tends to $(q''\circ h)h'^2+(q'\circ h)h''=(q\circ h)''$ in
the norm of $C[a,b]$, as $n\to\infty$. There is an $n_0\in\NN$,
such that for all $n>n_0=n_0(\delta)$ we have $\norm{(p_{n}\circ
h)''-(q\circ h)''}{C[a,b]}<\delta$, but $(q\circ
h)''(x)=p_\delta''(x)\ge\delta$, $x\in[a,b]$, hence, $(p_{n}\circ
h)''(x)\ge0$, $x\in[a,b]$, and $p_{n}\circ h\in\Delta$. Clearly,
for all $n\in\NN$ $p_n\circ h\in A$, so, taking into account~(5)
for $j=0$, we get $q\circ h=p_\delta\in\overline{\Delta\cap A}$.
Taking $\delta$ arbitrary small, we obtain
$p\in\overline{\Delta\cap A}$.
\eop

Roughly speaking, it is sufficient to have at least one ``nice''
function in algebra to have convex approximation to all convex
functions on the segment. Moreover, Proposition~3 does not require
any shape of $h$, except for monotonicity.

Now we give a notation for the algebra of exponential polynomials.
Denote by $E_{d}$ the set of all functions $E:\RR^d\to\RR$
$$
E(\bfm{x})=\sum_{j=1}^{n}{c_j}e^{\bfm{\alpha_j\cdot x}}=
\sum_{j=1}^{n}{c_j}\exp\left\{\sum_{k=1}^d\alpha_{j,k} x_k\right\},
\quad \bfm x=(x_1,\ldots,x_d)\in\RR^d,
$$
for some $n\in\NN$, $c_j\in\RR$ and
$\bfm{\alpha_j}=(\alpha_{j,1},\ldots,\alpha_{j,d})\in\ZZ_+^d$,
$j=1,\ldots,n$, where $\ZZ_+$ is the set of all non-negative
integers.

In particular, the previous result shows, that for any
$[a,b]\subset\RR$, we have $\Delta=\overline{E_1\cap\Delta}$,
where $\Delta$ is the set of all continuous and convex functions
on $[a,b]$. Theorem~1a) allows to extend this to multivariate
approximation, namely, we have
\proclaim Corollary.  Let $K\subset \RR^d$ be a convex compact
set, and $\Delta$ be the set of convex and continuous functions on
$K$. Then $\Delta=\overline{E_d\cap\Delta}$.

\pf There is a sufficiently large $N$ such that
$K\subset[-N,N]^{d}$. Take $h(x)=e^x$, $x\in[-N,N]$. By the
previous, we have that functions $l(x):=x$, $x\in[-N,N]$ and $-l$
are in the closure of $E_1\cap\Delta_N$, where $\Delta_N$ is the
set of continuous and convex functions on $[-N,N]$. Thus, for all
$j=1,\ldots,d+1$ the function $l_j$ (see the formulation of
Theorem~1) is in the closure of $E_d\cap\Delta$. By Theorem~1a),
$\Delta=\overline{E_d\cap\Delta}$.
\eop

So, as an example of application of Theorem~1, we have showed that
any convex and continuous multivariate function on a convex
compact can be arbitrary well approximated by multivariate
exponential polynomials, which are also convex on the compact.


\Acknowledgments{Both authors thank the referees for their
constructive remarks, which led to noticeable improvement of the
paper. The first author is supported by Ukrainian-French grant
``Dnipro'', M/262-2003. The second author thanks the organizers of
``Advances in Constructive Approximation'' for the support for the
attending the conference.}


\References

\ref DeVore, R. A., G. G. Lorentz, {\sl Constructive
Approximation}, Springer Verlag, Berlin, 1993.

\ref Lorentz, G. G., M.\ v.\ Golitschek, Y.\ Makovoz,
{\sl Constructive Approximation}, Springer Verlag, Berlin, 1996.

\ref Shvedov, A. S., Coconvex approximation of functions in many
variables, Mat.~Sbornik, {\bf 115}, no.~4 (1981), 577--589.


{

\bigskip\obeylines
Andriy V. Bondarenko
Dept. of Math. Analysis
Faculty of Mech. and Mathematics
Kyiv National Taras Shevchenko University
Kyiv, 01033, Ukraine

{\tt andriybond@gmail.com}

\bigskip
Andriy V. Prymak
Department of Mathematics
University of Manitoba
Winnipeg, R3T2N2, Canada

{\tt prymak@gmail.com
http://prymak.net/}

}  

\bye